\documentclass[11pt]{article}
\usepackage{amsfonts}

\usepackage{amssymb}
\usepackage{latexsym}
\usepackage{amscd}
\usepackage[centertags]{amsmath}
\usepackage{xypic}
\newtheorem{defn}{Definition}[section]
\newtheorem{thm}{Theorem}
\newtheorem{cor}[thm]{Corollary}
\newtheorem{lem}[thm]{Lemma}
\newtheorem{prop}[thm]{Proposition}
\newtheorem{rem}{Remark}

\numberwithin{equation}{section}

\newcommand{\calF}{\mathcal{F}}

\newcommand{\rmd}{\mathrm{d}}

\newcommand{\Enum}{\mathbb{E}}
\newcommand{\Pnum}{\mathbb{P}}
\newcommand{\Qnum}{\mathbb{Q}}
\newcommand{\Rnum}{\mathbb{R}}

\newcommand{\Xnum}{\mathbb{X}}
\newcommand{\Znum}{\mathbb{Z}}
\newcommand{\Nnum}{\mathbb{N}}
\newcommand{\pr}{\mbox{\boldmath$p$}}

\newcommand{\Om}{\Omega}
\newcommand{\om}{\omega}
\newcommand{\lam}{\lambda}
\newcommand{\Lam}{\Lambda}
\newcommand{\vep}{\varepsilon}
\newcommand{\qed}{\hfill $\Box$}

\begin{document}
\title{Range-Renewal Structure in Continued Fractions} 

\author{Jun Wu$^{1}$ and Jian-Sheng Xie$^{2}$\footnote{Corresponding author.}\\
{\footnotesize 1. Department of Mathematics, Huazhong University of
Science and Technology, Wuhan 430074, P. R. China}\\
{\footnotesize 2. School of Mathematical Sciences, Fudan
University, Shanghai 200433, P. R. China}\\
{\footnotesize E-mails: jun.wu@mail.hust.edu.cn, jiansheng.xie@gmail.com}} 
\date{}
\maketitle
\begin{abstract}
Let $\om=[a_1, a_2, \cdots]$ be the infinite expansion of continued fraction for an irrational
number $\om \in (0,1)$; let $R_n (\om)$ (resp. $R_{n, \, k} (\om)$, $R_{n, \, k+} (\om)$) be
the number of distinct partial quotients each of which appears at least once (resp. exactly $k$
times, at least $k$ times) in the sequence $a_1, \cdots, a_n$. In this paper it is proved that
for Lebesgue almost all $\om \in (0,1)$ and all $k \geq 1$,
$$
\displaystyle \lim_{n \to \infty} \frac{R_n
(\om)}{\sqrt{n}}=\sqrt{\frac{\pi}{\log 2}}, \quad \lim_{n \to
\infty} \frac{R_{n, \, k} (\om)}{R_n (\om)}=\frac{C_{2 k}^k}{(2k-1)
\cdot 4^k}, \quad \lim_{n \to \infty} \frac{R_{n, \, k} (\om)}{R_{n,
\, k+} (\om)}=\frac{1}{2k}.
$$
The Hausdorff dimensions of certain level sets about $R_n$ are
discussed.
\end{abstract}

\section{Introduction}\label{sec:1}
In early 2011, a beautiful range-renewal structure in i.i.d. models
was found by Chen et al \cite{CXY}. Among other results, the typical
main results in \cite{CXY} say that, given $n$ samples of a heavy
tailed \textbf{regular} (see \cite{CXY} for the definition) discrete distribution $\pi$ with an
intrinsic index $\gamma=\gamma (\pi) \in (0,1)$, then
\begin{eqnarray}
\displaystyle \lim_{n \to \infty} \frac{R_n}{\Enum R_n} &=& 1,\\
\lim_{n \to \infty} \frac{R_{n, \, k}}{R_n} &=& r_{_k} (\gamma)
:=\frac{\gamma \cdot \Gamma (k-\gamma)}{k! \cdot \Gamma (1-\gamma)}, \\
\lim_{n \to \infty} \frac{R_{n, \, k}}{R_{n, \, k+}} &=&
\frac{\gamma}{k},
\end{eqnarray}
where $R_n$ (resp. $R_{n, \, k}$, $R_{n, \, k+}$) stands for the number of distinct sample values
each of which appears at least once (resp. exactly $k$ times, at least $k$ times). And the so-called
range-renewal speed $\Enum R_n$ can be calculated out explicitly in $n$; for
instance, if
$$
\pi_x=\frac{C}{x^\alpha} \cdot [1+o(1)], x \in \Nnum
$$
with $1<\alpha<+\infty$, then
$$
\gamma=\gamma (\pi) :=1/\alpha \hbox{ and } \Enum R_n=\Gamma (1-\gamma) \cdot (C
\cdot n)^\gamma \cdot [1+o (1)].
$$

Soon after that, it was realized by the second author that the above
results should be somehow universal. Especially it should be true in
the continued fractions system (equipped with the well-known Gauss
measure $\mu$) since this system, in a certain sense, is stationary and
strongly mixing, meaning that the system is very near an i.i.d. model.
The confidence was strengthened in the spring of 2012 after a numerical
simulation with the help of Mr. Peng Liu.

Given an irrational number $\om \in (0,1)$. Let $R_n (\om)$ (resp.
$R_{n, \, k} (\om)$, $R_{n, \, k+} (\om)$) be the number of distinct
partial quotients each of which appears at least once (resp. exactly
$k$ times, at least $k$ times) in the first $n$ partial
quotients of $\om$. In this paper, we shall prove the following
interesting result: for Lebesgue almost all $\om \in (0,1)$ and all $k \geq 1$
$$
\frac{R_n (\om)}{\sqrt{n}}=\sqrt{\frac{\pi}{\log 2}}+ o(1),
\frac{R_{n, \, k} (\om)}{R_n (\om)}=\frac{C_{2 k}^k}{(2k-1)
\cdot 4^k}+o (1)  \hbox{ and } \frac{R_{n, \, k} (\om)}{R_{n, \, k+}
(\om)}=\frac{1}{2k}+o (1)
$$
as $n \to +\infty$; see Theorem \ref{thm: main} for the explicit
statement and its proof in Section \ref{sec:3}. Moreover we would
discuss the Hausdorff dimension of certain level sets; see Theorem
\ref{thm: main2} and its proof in Section \ref{sec:4}. As is pointed
in a remark (Remark \ref{rem: 1}), though the continued fraction
system (equipped with the Gauss measure $\mu$) is in fact a positive
recurrent system, we can observe certain kind of escape phenomenon:
For any $k \in \Nnum$ and Lebesgue almost all $\om \in (0, 1)$
$$
\lim_{n \to +\infty} \frac{R_{n, \, k} (\om)}{R_{n, \, k+}
(\om)}=\frac{1}{2k}.
$$
In the simple symmetrical random walk model in $\Znum^d$ (with $d
\geq 3$), the above limited ratio is always the escape rate
$\gamma_d$ \cite{DE50, ET, Revesz} where $R_n$ (respectively, $R_{n, \, k},
R_{n, \, k+}$) is interpreted as the number of distinct cites
visited at least once (respectively, exactly $k$ times, at least $k$
times) up to time $n$.

It is worthwhile to point out that, the ideas in this paper are
applicable to other systems to obtain similar results to Theorem
\ref{thm: main}. In view of the techniques developed in this article,
it's also possible to obtain further results (as those in \cite{CXY}) for the
current continued fractions model.
\section{Main Settings and Results}\label{sec:2}
Throughout this paper, the notation $y=O (z)$ implies that there
exists some universal constant $C>0$ such that $C^{-1} \leq
|\frac{y}{z}| \leq C$; the notation $y=\overline{O} (z)$ implies
that there exists some universal constant $C>0$ such that
$|\frac{y}{z}| \leq C$. And the notation $y=o (z)$ is understood in
the usual way. We shall use $O (1)$ to denote universal constants
which could change from line to line. For two sets $A, B$, we will
write $AB :=A \bigcap B$ for simplicity.

Let $\Xnum=(0,1) \setminus \Qnum$ be the set of irrational numbers in
the interval $(0,1)$. For any $\om \in \Xnum$, let $\{a_n=a_n
(\om)\}_{n=1}^\infty$ be the partial quotients of $\om$ in continued
fraction form, i.e.,
$$
\om=[a_1, a_2, \cdots]
:=\frac{1}{a_1+\frac{1}{a_2+\frac{1}{\ddots}}}.
$$
Therefore for any $\om \in \Xnum$ there is a unique (natural) coding
$(a_1, a_2, \cdots)$ (still written $\om$ for simplicity) in
$\Om=\Nnum^\Nnum$. And the Gauss map $T: \Xnum \to \Xnum$
$$
T (\om) :=\frac{1}{\om} \; (\mathrm{mod\;} 1)=[a_2, a_3, \cdots]
$$
induces the natural left-shift map $\sigma: \Om \to \Om$. The Gauss
measure $\mu$ (which satisfies $\rmd \mu (\om)=\frac{\rmd \om}{(\log
2) \cdot (1+\om)}$ and which is invariant under $T$) on $(0,1)$
naturally induces a probability measure $\Pnum$ on $\Om$. And for
any $x \in \Nnum$ we know
\begin{equation}\label{eq: def4pi}
\pi_x :=\Pnum (a_1=x)=-\log_2
[1-\frac{1}{(x+1)^2}]=\frac{1}{(\log 2) \cdot (x+1)^2}+O
(\frac{1}{x^4})
\end{equation}
as $x \to +\infty$. Therefore there is a probability measure
$\pi=(\pi_x: x \in \Nnum)$ on $\Nnum$ which also naturally induces
an infinitely independent product measure $\widetilde{\Pnum}
:=\pi^\infty$ on $\Om$. The expectation operator of the probability
measure $\Pnum$ (respectively, $\widetilde{\Pnum}$) will be denoted
by $\Enum$ (respectively, $\widetilde{\Enum}$). And we have the
following commuting graph
$$
\begin{array}{rcl}
(\Om, \Pnum) &\stackrel{\sigma}{\longrightarrow}& (\Om, \Pnum)\\
\pr \downarrow \; & & \; \downarrow \pr\\
(\Xnum, \mu) &\stackrel{T}{\longrightarrow}& (\Xnum, \mu)
\end{array}
$$
with $\pr$ being the natural projection
$$
\pr (a_1, a_2, \cdots) :=[a_1, a_2, \cdots].
$$
Due to this obvious identification, we shall \textit{not} distinguish the
spaces $\Om$ and $\Xnum=(0,1)\setminus \Qnum$ from hereon.

Given $\om \in \Xnum$. For any fixed $x \in \Nnum$, write
\begin{equation}\label{eq:def4N-n}
N_n (x)=N_n (x, \, \om) :=\sum_{k=1}^n 1_{\{a_k (\om)=x\}}
\end{equation}
which is the visiting number of the state $x$ by the partial quotients $a_1 (\om), \cdots, a_n (\om)$. Define
\begin{equation}\label{eq:def4R-n}
R_n (\om) :=\sum_{x \in \Nnum} 1_{\{N_n (x, \, \om) \geq 1\}};
\end{equation}
this is the number of distinct values of $a_1 (\om), \cdots, a_n (\om)$, i.e., $R_n (\om)=\sharp \{a_1 (\om), \cdots, a_n (\om)\}$.
Define also for any $k \in \Nnum$
\begin{equation}\label{eq:def4R-n-k}
R_{n, \, k} (\om) :=\sum_{x \in \Nnum} 1_{\{N_n (x, \, \om) =k\}};
\end{equation}
this is the number of distinct partial quotients each of which appears
exactly $k$ times in the finite sequence $a_1 (\om), \cdots, a_n
(\om)$.

Our main result is as the following.
\begin{thm}\label{thm: main}
For Lebesgue almost all $\om \in (0,1)$, we have
\begin{equation}
\lim_{n \to \infty} \frac{R_n (\om)}{\sqrt{n}}=\sqrt{\frac{\pi}{\log
2}}.
\end{equation}
Furthermore, for any $k \geq 1$
\begin{equation}
\lim_{n \to \infty} \frac{R_{n, \, k} (\om)}{R_n (\om)}=\frac{C_{2
k}^k}{(2k-1) \cdot 4^k} =:r_{_k},
\end{equation}
where $C_n^m:=\frac{n!}{m! (n-m)!}$ and $\sum \limits_{k=1}^\infty r_{_k}=1$.
\end{thm}

\begin{rem}\label{rem: 1}
From the above theorem, we know
\begin{itemize}
  \item[{\rm (1)}] as $k \to \infty$, $r_{_k}=\frac{1}{2\sqrt{\pi} } \cdot
  k^{-3/2}+O (k^{-5/2})$, which is a power law with index $3/2$;
  \item[{\rm (2)}] Let
\begin{equation}\label{eq:def4R-n-k+}
R_{n, \, k+} (\om) :=\sum_{\ell=k}^n R_{n, \, \ell} (\om).
\end{equation}
Then for Lebesgue almost all $\om \in (0,1)$
\begin{equation}\label{eq: escape-rate}
\lim_{n \to +\infty} \frac{R_{n, \, k} (\om)}{R_{n, \, k+}
(\om)}=\frac{1}{2k}
\end{equation}
which can be interpreted as an average escape rate at the level $k
\geq 1$, though the model itself is in fact positive recurrent. This
result can be seen as the following. Put for $k \geq 2$
$$
r_{_{k+}}=2k \cdot r_{_k}=\frac{2k}{2k-1} \cdot
\frac{C_{2k}^k}{4^k}=\prod_{j=1}^{k-1} (1-\frac{1}{2 j })
$$
and $r_{_{1+}}=1$. Obviously
$$
r_{_{k+}}-r_{_{(k+1)+}}=r_{_{k+}}-r_{_{k+}} \cdot
(1-\frac{1}{2k})=\frac{1}{2k} \cdot r_{_{k+}}=r_{_k}, \quad k
\geq 1.
$$
Thus $\displaystyle r_{_{k+}}=\sum_{\ell \geq k} r_{_\ell} $. Then
the result follows from Theorem \ref{thm: main}.

  \item[{\rm (3)}] We could recall a standard result in the field
of the simple symmetrical random walk (abbr. \textbf{SSRW}) on
$\Znum^d$ (with $d \geq 3$), though it was not expressed explicitly
in \cite{ET} (but nearly explicitly as Theorem 20.11 in \cite[p. 220]{Revesz}).
It says that the limited ratio in (\ref{eq: escape-rate}) is always $\gamma_d$,
the usual escape rate of the random walk. The proof of this result has only
several lines using sub-additive ergodic theorem \cite{Kingman68, Kingman73, Kingman76}
(essentially the argument of Derriennic \cite{Derriennic}) which we
would like to list as the following. In \textbf{SSRW} on $\Znum^d$
(with $d \geq 3$), $R_n, \,R_{n, \, k-} :=R_n-R_{n, (k+1)+}$ are all
sub-additive. Hence
\begin{eqnarray*}
\lim_{n \to +\infty} \frac{R_n}{n} &=& \lim_{n \to +\infty}
\frac{\Enum R_n}{n}=\gamma_d,\\
\lim_{n \to +\infty} \frac{R_{n, \, k-}}{n} &=& \lim_{n \to +\infty}
\frac{\Enum R_{n, \, k-}}{n} =\gamma_d \cdot[1- (1-\gamma_d)^{k}].
\end{eqnarray*}
(Cf. \cite{XX} for more detailed calculations.) And the result follows.
This is a different approach for Dvoretzky and Erd\"{o}s' result \cite{DE50}
for $d \geq 3$; The variation estimations are not needed in this proof.
\end{itemize}
\end{rem}

\begin{rem}\label{rem: 1'}
In view of Corollary \ref{cor: main-tool}, the result in the above theorem can be strengthened
as the following: for any $k \geq 1$ fixed
\begin{eqnarray}
\frac{R_n (\om)}{\sqrt{n}} &=& \sqrt{\frac{\pi}{\log 2}}+ o (\frac{1}{n^{0.0698}}),\\
\frac{R_{n, \, k} (\om)}{R_n (\om)} &=& \frac{C_{2 k}^k}{(2k-1)
\cdot 4^k}+o (\frac{1}{n^{0.0698}}),  \\
\frac{R_{n, \, k} (\om)}{R_{n, \, k+} (\om)} &=& \frac{1}{2k}+o (\frac{1}{n^{0.0698}})
\end{eqnarray}
almost surely since $\frac{\vep_0}{3}=\frac{1}{6 (1+2 \log 2)}=0.069843...$.
\end{rem}

Put for any $\beta \geq 0, c>0$
$$
E (\beta, c) :=\{\om \in (0,1) \setminus \Qnum: \lim_{n \to +\infty}
\frac{R_n (\om)}{n^\beta}=c\}.
$$
The above theorem implies that
$$
\dim_H E (\beta, c)=1 \hbox{ for } \beta=\frac{1}{2} \hbox{ and }
c=\sqrt{\frac{\pi}{\log 2}}.
$$
One would wonder what happens for other choices of $(\beta, c)$.

Let
$$
d_H (\beta):=\dim_H \{\om \in (0,1) \setminus \Qnum: 0<\lim_{n \to
+\infty} \frac{R_n (\om)}{n^\beta}<+\infty\}.
$$
The above theorem proves $d_H (\frac{1}{2})=1$. The classical result
of Jarn\'{\i}k (see \cite{Ja28}) implies that $d_H (0)=1$. Therefore a natural
conjecture seems to be
$$
d_H (\beta)=1, \forall \beta \in [0,1).
$$
In fact we have the following result.
\begin{thm}\label{thm: main2}
Let $E (\beta, c)$ be defined as above and put $F (c)=E (1, c)$,
i.e.,
\begin{eqnarray*}
E (\beta, c)&:=& \Big\{ \om \in [0,1):\lim_{n\to
\infty}\frac{R_n(\om)}{n^\beta}=c \Big\},\\
F (c) &:=& \Big\{ \om \in [0,1):\lim_{n\to
\infty}\frac{R_n(\om)}{n}=c \Big\}.
\end{eqnarray*}
Then
\begin{itemize}
  \item[{\rm (1)}] for any $\beta \in (0, 1)$ and $c>0$
\begin{equation}\label{eq: dim-1}
\dim_H E (\beta, c)=1;
\end{equation}
  \item[{\rm (2)}] for any $c \in (0, 1]$
\begin{equation}\label{eq: dim-2}
\dim_H F (c)=\frac{1}{2};
\end{equation}
\end{itemize}
\end{thm}
\begin{rem}\label{rem:2}
Actually, the proof of Theorem \ref{thm: main2} in Section
\ref{sec:4} can be modified to prove the following more general
result. For any smooth function $\psi$ with $\psi (n) \nearrow \infty,
\frac{\psi (n)}{n} \searrow 0$ as $n \to \infty$, the set
$$
E_{\psi} :=\{\om \in [0,1): \lim_{n \to +\infty} \frac{R_n
(\om)}{\psi (n)}=1\}
$$
is always of Hausdorff dimension $1$.

It is interesting to point out that $e-2=0.71828... \in F (\frac{1}{3})$, which is due to Euler (Cf.
\cite[pp. 12]{I-K}); one might expect $\pi -3=0.14159...$ satisfying the conclusions of Theorem \ref{thm: main}.
\end{rem}
The above results in Theorem \ref{thm: main2} and in Remark \ref{rem:2} fall into the so-called
\textbf{fractional dimensional theory}. This theory is attained much attention in studying the exceptional
sets arising in the metrical theory of continued fractions. It seems that the first published work
in this region is a paper by Jarn\'{i}k \cite{Ja28}. Later on, Good \cite{Go} gave an quite
overall investigation of sets with some restrictions on their partial quotients. For more results
among this area, one can refer to the work of Hirst \cite{Hi}, L\'uczak \cite{Lu}, Mauldin and
Urb\'{a}nski \cite{MU1}, \cite{MU2}, Pollicott and Weiss \cite{PW}, Wang and Wu \cite{WW},  Li,
Wang, Wu and Xu \cite{LWWX} and references therein.
\section{Proof of Theorem \ref{thm: main}}\label{sec:3}
The proof somehow follows the main strategy of \cite{DE50}. The main
tool in use is the one developed in \cite{CXY} (also traced back to \cite{DE50} in some sense)
which we restate here without proof as the following.
\begin{lem}\label{lem: main-tool}
Let $\{Y_n\}_1^\infty$ be a sequence of non-negative random
variables in a probability space $(\Om, \Pnum)$. Let $\displaystyle
S_n :=\sum_{k=1}^n Y_n$. Suppose $\displaystyle \lim_{n \to +\infty}
\Enum S_n=+\infty, \sup\{\Enum Y_n: n \geq 1 \}<+\infty$ and
\begin{equation}\label{eq: cond4SLLN}
\mathrm{Var\,} (S_n) \leq C \cdot (\Enum S_n)^{2-\delta}
\end{equation}
for some $\delta>0, C>0$ and all sufficiently large $n$. Then
\begin{equation}\label{eq: SLLN}
\lim_{n \to \infty} \frac{S_n}{\Enum S_n}=1
\end{equation}
holds true almost surely. The condition (\ref{eq: cond4SLLN}) can
even be weakened as
\begin{equation}\label{eq: cond4SLLN'}
\mathrm{Var\,} (S_n) \leq C \cdot (\Enum S_n)^2/(\log \Enum
S_n)^{1+\delta}.
\end{equation}
\end{lem}

\begin{cor}\label{cor: main-tool}
The conclusion of the above lemma can be strengthened as the following: for any fixed $\beta
\in (0, \frac{\delta}{3})$
\begin{equation}\label{eq: SLLN-1}
\frac{S_n}{\Enum S_n}=1+\bar{O} (\frac{1}{(\Enum S_n)^{\beta}})
\end{equation}
holds true almost surely if condition (\ref{eq: cond4SLLN}) holds. If condition
(\ref{eq: cond4SLLN'}) holds instead of (\ref{eq: cond4SLLN}), then
\begin{equation}\label{eq: SLLN-1}
\frac{S_n}{\Enum S_n}=1+\bar{O} (\frac{1}{(\log \Enum S_n)^{\beta}}).
\end{equation}
\end{cor}

For the convenience of the readers, we translate the main results in
\cite{CXY} for the probability measure $\widetilde{\Pnum}$ induced
by i.i.d. distribution $(\pi_x, x \in \Nnum)$ (see (\ref{eq:
def4pi}) for the definition) as below, which is in fact also the
motivation of the current paper.
\begin{lem}\label{lem: basic-facts}
We have the following estimations for the measure
$\widetilde{\Pnum}$ as $n \to +\infty$.
\begin{itemize}
  \item[{\rm (1)}] $\displaystyle \widetilde{\Enum} R_n=\sum_{x=1}^\infty
\Bigl[ 1-(1-\pi_x)^n \Bigr]=\sqrt{\frac{\pi n}{\log 2}}
+\overline{O} (1), \quad \mathrm{Var}_{\widetilde{\Pnum}} (R_n) \leq
\widetilde{\Enum} R_n$;
  \item[{\rm (2)}] For any fixed $k \geq 1$,
$$
\widetilde{\Enum} R_{n, \, k}=\sum_{x=1}^\infty C_n^k \cdot \pi_x^k
\cdot (1-\pi_x)^{n-k}=\sqrt{\frac{\pi n}{\log 2}} \cdot r_{_k}
+\overline{O} (1), \quad \mathrm{Var}_{\widetilde{\Pnum}} (R_{n, \,
k}) \leq [1+o (1)] \cdot \widetilde{\Enum} R_{n, \, k},
$$
where $r_{_k}$ is defined in Theorem \ref{thm: main}. Let
$\displaystyle r_{_{k+}} := \sum_{\ell=k}^\infty
r_{_\ell}=1-\sum_{\ell=1}^{k-1} r_{_\ell}$, then
$$
\widetilde{\Enum} R_{n, \, k+}=\sqrt{\frac{\pi n}{\log 2}} \cdot
r_{_{k+}} +\overline{O} (1);
$$
  \item[{\rm (3)}] For $\widetilde{\Pnum}$-a.e. $\om$ and all $k \geq 1$,
$\displaystyle \lim_{n \to +\infty} \frac{R_n (\om)}{\widetilde{\Enum} R_n}=1,
\quad \lim_{n \to +\infty} \frac{R_{n, \, k} (\om)}{R_n
(\om)}=r_{_k}$.
\end{itemize}
\end{lem}

\subsection{Estimating $\Enum R_n$}\label{sec:3.1}
In this part we shall estimate $\Enum R_n$. By the definition
(\ref{eq:def4R-n}), we need to estimate the probability
$$
\Pnum (N_n (x) \geq 1)=\Pnum (a_k (\om)=x \hbox{ for some } 1 \leq k
\leq n)
$$
for any $x \in \Nnum$. Actually we will do this for large enough $x$
at first, i.e., for those
\begin{equation}
x > x^*_n :=\left\lfloor \sqrt{\frac{C n}{\log n}} \right \rfloor
\end{equation}
with $C>0$ to be determined later; here $\left\lfloor a \right \rfloor$
denotes the integer part of a real number $a$.

The basic idea is to exploit the following standard result in probability
theory.

\begin{lem}\label{lem: Basic-1}
Let $\{A_k\}_1^n$ be a sequence of measurable sets in a probability
space $(\Om, \calF, \Pnum)$. Then
\begin{equation}\label{eq: basic-eq}
\Pnum (\bigcup_{k=1}^n A_k)=\sum_{k=1}^n (-1)^{k-1} \sum_{1 \leq
i_1<\cdots<i_k \leq n} \Pnum (\bigcap_{r=1}^k A_{i_r}).
\end{equation}
\end{lem}
When taking the probability $\Pnum$ as a Dirac point measure, the
above equation is just an identity for indicator functions:
\begin{equation}\label{eq: basic-eq'}
1_{\{N_A \geq 1\}}=\sum_{k=1}^n (-1)^{k-1} \sum_{1 \leq
i_1<\cdots<i_k \leq n} 1_{A_{i_1, \cdots, i_k}}.
\end{equation}
From hereon, we will always write
\begin{equation}\label{eq: def4NA}
N_A =N_A (\om) :=\sum_{k=1}^n 1_{A_k} (\om), \quad A_{i_1, \cdots, i_k}
:=\bigcap_{r=1}^k A_{i_r}
\end{equation}
for a sequence $A=\{A_k\}_1^n$ of measurable sets.

Obviously, in order to estimate $\Pnum (N_n (x) \geq 1)$ by exploiting the formula (\ref{eq: basic-eq}),
one has to take in count the strong mixing property of the model under investigation. Therefore we would
introduce this property in the sequel.

For any $n\geq 1$ and $(a_1,\cdots, a_n)\in \Nnum^n$, call
\begin{eqnarray*}
I(a_1,\cdots, a_n)= \left\{
  \begin{array}{ll}
    \displaystyle \bigg[\frac{p_n}{q_n}, \frac{p_n+p_{n-1}}{q_n+q_{n-1}}\bigg) ,
& \hbox{{\rm{when}} \ $n$ {\rm{is \ even}}} \\
   \displaystyle \bigg(\frac{p_n+p_{n-1}}{q_n+q_{n-1}}, \frac{p_n}{q_n}\bigg],
& \hbox{{\rm{when}} \ $n$ {\rm{is \ odd}}}
  \end{array}
\right.
\end{eqnarray*}
an $n$-th order cylinder, where $\{p_{_k}, \ q_{_k}\}_{k=1}^n$ are
determined by following recursive relations
\begin{eqnarray}\label{f1}
p_{_k}=a_{_k} \cdot p_{_{k-1}}+p_{_{k-2}}, \  q_{_k}= a_{_k} \cdot q_{_{k-1}}+q_{_{k-2}},
\ 1\leq k \leq n
\end{eqnarray}
with the conventions that $p_{_{-1}}=1, p_{_0}=0$, $q_{_{-1}}=0, q_{_0}=1$. It
is well known, see \cite{Kh63}, that $I(a_1,\cdots, a_n)$ just
represents the set of points in $[0,1)$ which have a continued
fraction expansions beginning with $a_1,\cdots, a_n$.
\begin{prop}\label{p1} (See \cite{Kh63}.)
For any $n\geq 1$ and $(a_1,\cdots,a_n)\in \mathbb{N}^n$, one has
\begin{eqnarray}\label{f2}
\big|I(a_1,\cdots, a_n)\big|=\frac{1}{q_n(q_n+q_{n-1})},
\end{eqnarray}
where $|I(a_1,\cdots, a_n)|$ denotes the length of $I(a_1,\cdots,
a_n)$.
\end{prop}

The following lemma is a standard result in continued fraction
model, which is our desired mixing property.
\begin{prop}\label{prop: mixing} (See \cite{Billingsley} or \cite{I-K}.)
In the continued fraction model,
\begin{equation}\label{eq: mixing-eq}
|\frac{\mu (I (\tilde{x}) \bigcap T^{-(m+L)} I (\tilde{y}))}{\mu (I (\tilde{x})) \cdot
\mu (I (\tilde{y}))} -1| \leq O (q^L)
\end{equation}
for some $q \in(0, 1)$, all (fixed) $\tilde{x}=(x_1, \cdots, x_m), \tilde{y}=(y_1,
\cdots, y_{n})$ and sufficiently  large $L$.
\end{prop}
The first contributor for the mixing property of continued fractions
is Kuzmin \cite{Kuzmin} (see also \cite{Kh63, Kh64}), who proved a
sub-exponential decay rate in solving Gauss's conjecture on
continued fractions. L\'{e}vy \cite{Levy29} (see also \cite[Ch.
IX]{Levy54}) independently proved the exponential decay rate with $q=3.5-2
\sqrt{2} =0.67157...$, also solving Gauss's conjecture. Using Kuzmin's approach,
Sz\H{u}sz \cite{Szusz} claimed to have lowered the L\'{e}vy estimate for $q$
to $0.4$ whereas his argument yields $q=0.485$ rather than $q=0.4$.
The optimal value of $q =0.30366300289873265859...$ was determined
by Wirsing \cite{Wirsing}.

The above equations alone are not sufficient. We need the
following observation.
\begin{lem}\label{lem: continued-fraction-1}
For any $m, n \geq 1$ and $\tilde{x}=(x_1, \cdots, x_m) \in \Nnum^m,
\tilde{y}=(y_1, \cdots, y_n) \in \Nnum^n$
\begin{equation}\label{eq: propability-bound}
K^\prime :=\log2 \leq \frac{\mu (I (\tilde{x}, \tilde{y}))}{\mu (I(\tilde{x})) \cdot \mu (I(\tilde{y}))}
\leq K:=2 \log 2,
\end{equation}
where the above bounds are optimal.
\end{lem}
We can compare (\ref{eq: propability-bound}) with a standard result listed below
\begin{equation}\label{eq: propability-bound'}
\frac{1}{2} \leq \frac{|I (\tilde{x}, \tilde{y})|}{|I(\tilde{x})| \cdot |I(\tilde{y})|} \leq 2,
\end{equation}
both bounds of which are also optimal; this can be easily proved in view of Proposition \ref{p1}.

\noindent \textit{Proof of Lemma \ref{lem: continued-fraction-1}}. \;
Let
$$
\frac{p_k}{q_k}=[x_1, \cdots, x_k], k=1, \cdots, m
$$
be irreducible fractions. Then we clearly have
$$
\frac{|I (\tilde{x})|}{\log 2}  \cdot \min_{\om \in I (\tilde{x})} \frac{1}{1+\om}\leq \mu (I (\tilde{x}))=\int_{I (\tilde{x})} \frac{\rmd
\om}{(\log 2) \cdot (1+\om)} \leq \frac{|I (\tilde{x})|}{\log 2} \cdot \max_{\om \in I (\tilde{x})} \frac{1}{1+\om}
$$
and $|I (\tilde{x})|=\frac{1}{q_m \cdot (q_m+q_{m-1})}$ by Proposition \ref{p1}. Furthermore,
\begin{eqnarray*}
\mu (I (\tilde{x}, \tilde{y})) &=& \int_{I (\tilde{x}, \tilde{y})} \frac{\rmd \om}{(\log 2) (1+\om)}\\
&=& \int_{I (\tilde{y})} \frac{\rmd \om}{(\log 2) (1+\om)} \cdot \frac{1+\om}{(q_m+\om \cdot q_{m-1})
\cdot [q_m+p_m+\om \cdot (q_{m-1}+p_{m-1})]}\\
&\leq& \mu (I (\tilde{y})) \cdot \max_{\om \in I (\tilde{y})} \frac{1+\om}{(q_m+\om \cdot q_{m-1})
\cdot [q_m+p_m+\om \cdot (q_{m-1}+p_{m-1})]}.
\end{eqnarray*}
Therefore
$$
\frac{\mu (I (\tilde{x}, \tilde{y}))}{\mu (I (\tilde{x})) \cdot \mu (I (\tilde{y}))} \leq
K (\tilde{x}, \tilde{y}) \cdot \log 2,
$$
where
\begin{eqnarray*}
K (\tilde{x}, \tilde{y}) &:=& \max\{\frac{(1+\om) \cdot (1+\om^\prime) [ q_m
\cdot (q_m+q_{m-1})]}{(q_m+\om \cdot q_{m-1}) \cdot [q_m+p_m+\om \cdot (q_{m-1}+p_{m-1})]}:
\om^\prime \in I (\tilde{x}), \om \in I (\tilde{y}) \}\\
&=& \max_{\om \in I (\tilde{y})} \{\frac{(1+\om) (q_m+q_{m-1}) (p_m+q_m)}{(q_m+\om \cdot q_{m-1})
\cdot [q_m+p_m+\om \cdot (q_{m-1}+p_{m-1})]},
\\
&& \frac{(1+\om) \cdot q_m (p_m+q_m+p_{m-1}+q_{m-1})}{(q_m+\om \cdot q_{m-1}) \cdot [q_m+p_m+\om
\cdot (q_{m-1}+p_{m-1})]}\}.
\end{eqnarray*}
Obviously
$$
\frac{(1+\om) \cdot (q_m+q_{m-1})}{q_m+\om \cdot q_{m-1}} \leq 2 \hbox{ and } \frac{(1+\om)
(p_m+q_m+p_{m-1}+q_{m-1})}{q_m+p_m+\om \cdot (q_{m-1}+p_{m-1})} \leq 2,
$$
which admits $K=2 \cdot \log 2$; similarly one can take $K^\prime=\log 2$.
$K^\prime=\log 2, K=2 \log2$ is optimal since $\displaystyle
\lim_{a \to +\infty} \frac{\mu (I (a, a))}{\mu (I (a))^2}=\log 2$, $\displaystyle
\lim_{a \to +\infty} \frac{\mu (I (a, 1, a))}{\mu (I (a)) \cdot \mu (I (1, a))}=2 \log 2$.
\qed
\begin{rem}
Consider the optimal bounds of the following inequality
\begin{equation}\label{eq: propability-bound-2}
K^\prime_r \leq \frac{\mu (I (\tilde{x}_{_1}, \cdots, \tilde{x}_{_r}))}{\mu (I(\tilde{x}_{_1}))
\cdots \mu (I(\tilde{x}_{_r}))} \leq K_r,
\end{equation}
where $\tilde{x}_1, \cdots, \tilde{x}_r$ are all natural numbers tuples (possibly of different length)
and $r \geq 2$. A careful calculation reveals
$$
\lim_{a \to +\infty} \frac{\mu (I (\{a\}^r))}{\mu (I(a))^r}=(\log 2)^{r-1} \hbox{ and }
\lim_{a \to +\infty} \frac{\mu (I (\{1a\}^r))}{\mu (I(1a))^r}=(2 \log 2)^{r-1},
$$
where $\{a\}^r$ (resp. $\{1a\}^r$) denotes the integer tuple obtained by repeating the word ``$a$"
(resp. ``$1a$") $r$ times. Hence $K^\prime_r=(\log 2)^{r-1}$ and $K_r=(2 \log 2)^{r-1}$ are optimal.

Similarly, the following bounds are also optimal:
\begin{equation}\label{eq: propability-bound'}
\frac{1}{2^{r-1}} \leq \frac{|I (\tilde{x}_{_1}, \cdots, \tilde{x}_{_r})|}{|I(\tilde{x}_{_1})|
\cdots |I(\tilde{x}_{_r})|} \leq 2^{r-1}.
\end{equation}
\end{rem}

Now for a fixed integer $x >x^*_n$, we take $A_i=A^x_i=\{\om: a_i
(\om)=x\}, i=1, \cdots, n$ in eq. (\ref{eq: basic-eq}). Noting
$\Pnum (A^x_i)=\widetilde{\Pnum} (A^x_i)=\pi_x$, we would have
\begin{eqnarray*}
\Delta (x) &:=& |\Pnum (N_n (x) \geq 1)-\widetilde{\Pnum} (N_n (x)
\geq 1)| \\
&\leq& \sum_{k=2}^n \sum_{1 \leq i_1<\cdots<i_k \leq n} |\Pnum
(A^x_{i_1, \cdots, \, i_k})-\widetilde{\Pnum} (A^x_{i_1, \cdots, \,
i_k})|.
\end{eqnarray*}
Here $A^x_{i_1, \cdots, \, i_k}:=\{\om: a_{i_r} (\om)=x, r=1,
\cdots, k\}$.

We choose an integer number
$$
K_n := \left \lfloor C_1 \cdot \log n \right \rfloor
$$
for some sufficiently large $C_1$ (to be determined later). Consider
the sum
\begin{equation}
\Delta_1 (x) :=\sum_{k=K_n +1}^n \sum_{1 \leq i_1<\cdots<i_k\leq n}
|\Pnum (A^x_{i_1, \cdots, \, i_k}) -\widetilde{\Pnum} (A^x_{i_1,
\cdots, \, i_k})|.
\end{equation}
Clearly $\widetilde{\Pnum} (A^x_{i_1, \cdots, \, i_k})=\pi_x^k$ and in view of
Lemma \ref{lem: continued-fraction-1}
$$
\Delta_1 (x) \leq \sum_{k=K_n +1}^n C_n^k \cdot 2 \cdot K^{k-1} \cdot
\pi_x^k.
$$
Since $\pi_x \leq \frac{1}{\log 2} \cdot x^{-2}$, we have
\begin{eqnarray*}
\sum_{x >x^*_n} \Delta_1 (x) &\leq& \sum_{k=K_n +1}^n C_n^k \cdot 2
\cdot K^{k-1} \sum_{x >x^*_n} \pi_x^k\\
&\leq& \sum_{k=K_n +1}^n C_n^k \cdot 2 \cdot K^{k-1} \cdot O (1)
\cdot (\frac{1}{\log 2})^k \cdot \frac{1}{2k-1} \cdot (x^*_n)^{-(2k-1)}\\
&\leq& \sum_{k=K_n +1}^n \frac{n^k}{k!} \cdot (\frac{K}{\log 2})^{k}
\cdot O (1) \cdot \frac{1}{k} (\frac{\log n}{Cn})^{k-1/2}\\
&=& \sqrt{\frac{n}{\log n}} \cdot O (1) \cdot \sum_{k=K_n +1}^n
\frac{1}{(k+1)!} \cdot (\frac{K}{C \log 2} \cdot \log n)^k.
\end{eqnarray*}
Let $\alpha :=\frac{K}{C \log 2}=\frac{2}{C}$. By Stirling's formula
we see that if $C_1 \cdot C >4 e^2 \approx 29.55...$, then
\begin{eqnarray*}
\log (K_n!)-K_n \log (\alpha \log n) &=& \log \sqrt{2 \pi K_n} +K_n
\log K_n -K_n+o (1) -K_n \log (\alpha \log n)\\
&=& \log \sqrt{2 \pi K_n} +K_n \log \frac{K_n}{e \cdot \alpha \log
n}+o (1)\\
&\geq& C_1 \log n
\end{eqnarray*}
for sufficiently large $n$. And we will have
$$
\frac{\alpha \log n}{k}<\frac{1}{2}, \forall k>K_n.
$$
Therefore if $C_1 \geq 5$ and $C \geq 6$
\begin{eqnarray*}
\sum_{x >x^*_n} \Delta_1 (x) &\leq& \sqrt{\frac{n}{\log n}} \cdot O
(1) \cdot \sum_{\ell=1}^{n-K_n} \frac{(\alpha \cdot \log n)^{K_n}}{K_n!}
\cdot (\frac{1}{2})^\ell\\
&\leq& O (1) \cdot n^{-(C_1-\frac{1}{2})} \leq O (n^{-4}).
\end{eqnarray*}

Now choose a number $C_2 \geq 4$ and put
\begin{equation}
\lam_0 :=-\log q>0, \quad L_n :=\left \lfloor \frac{C_2}{\lam_0} \log n \right \rfloor+1.
\end{equation}
We want to estimate
\begin{equation}
\Delta_2 (x) :=\sum_{k=2}^{K_n} \sum_* |\Pnum (A^x_{i_1, \cdots, \,
i_k}) -\widetilde{\Pnum} (A^x_{i_1, \cdots, \, i_k})|
\end{equation}
and
\begin{equation}
\Delta_3 (x) :=\sum_{k=2}^{K_n} \sum_{**} |\Pnum (A^x_{i_1, \cdots,
\, i_k}) -\widetilde{\Pnum} (A^x_{i_1, \cdots, \, i_k})|,
\end{equation}
where the sum $\displaystyle \sum_*$ is over all $1 \leq
i_1<\cdots<i_k\leq n$ with $i_{u+1}-i_u >L_n, u=1, \cdots, k-1$, and
the sum $\displaystyle \sum_{**}$ is over the rest $1 \leq
i_1<\cdots<i_k\leq n$.

Clearly $q^{L_n} \leq n^{-C_2}$. Since $C_2 \geq 4$, for sufficiently large $n$
\begin{eqnarray*}
\Delta_2 (x) &\leq& \sum_{k=2}^{K_n} C_{n-(k-1)L_n}^k \cdot
\max[(1+O (q^{L_n}))^{k-1} -1, 1-(1-O (q^{L_n}))^{k-1}] \cdot \pi_x^k\\
&\leq& \sum_{k=2}^{K_n} \frac{n^k}{k!} \cdot O (n^{-3}) \cdot \pi_x^k.
\end{eqnarray*}
Thus if $C \geq 6>\frac{2}{\log 2}=2.88539...$, then
\begin{eqnarray*}
\sum_{x>x^*_n} \Delta_2 (x) &\leq& O (1) \cdot \sum_{k=2}^{K_n}
\frac{n^{k-3}}{k!} \cdot \sum_{x>x^*_n} \pi_x^k\\
&\leq& O (1) \cdot \sum_{k=2}^{K_n} \frac{n^{k-3}}{k!} \cdot
(\frac{1}{\log 2})^k \cdot \frac{1}{2k-1} \cdot (\frac{\log
n}{Cn})^{k-\frac{1}{2}}\\
&\leq& O (1) \cdot n^{-\frac{5}{2}} \sum_{k=2}^{K_n}
\frac{1}{(k+1)!} \cdot (\frac{\log n}{C \log 2})^k\\
&\leq& O (1) \cdot n^{-\frac{5}{2}} \cdot n^{\frac{1}{C \log 2}}
\leq O (n^{-2}).
\end{eqnarray*}

Now we estimate $\Delta_3 (x)$. Clearly
$$
C_n^k-C_{n-(k-1)L_n}^k \leq O (1) \cdot \frac{n^{k-1}}{k!} \cdot k^2
\cdot L_n.
$$
Therefore
$$
\Delta_3 (x) \leq \sum_{k=2}^{K_n} O (1) \cdot \frac{n^{k-1}}{k!}
\cdot k^2 \cdot L_n \cdot 2 K^{k-1} \cdot \pi_x^k.
$$
From hereon we take
\begin{equation}
C :=4+\frac{2}{\log 2}=6.88539...
\end{equation}
Then $\alpha :=\frac{K}{C \log 2}=\frac{\log 2}{1+2 \log 2}=0.290470...$.
For sufficiently large $n$
\begin{eqnarray*}
\sum_{x>x^*_n} \Delta_3 (x) &\leq& \sum_{k=2}^{K_n} O (1) \cdot
\frac{n^{k-1}}{k!} \cdot k^2 \cdot L_n \cdot (\frac{K}{\log 2})^k
\cdot \frac{1}{2k-1} \cdot (\frac{\log n}{Cn})^{k-\frac{1}{2}}\\
&\leq& O (1) \cdot \frac{L_n}{\sqrt{n \cdot \log n}} \cdot
\sum_{k=2}^{K_n} \frac{1}{(k-1)!} \cdot (\frac{K}{C \log 2} \cdot
\log n)^k\\
&\leq& O (1) \cdot \frac{(\log n)^{1.5}}{\sqrt{n}} \cdot n^{\alpha} \\
&=& O (\frac{(\log n)^{1.5}}{n^{\vep_0}}) \hbox{ with } \vep_0 :=\frac{1}{2 (1+2\log 2)}=0.209529....
\end{eqnarray*}

A similar calculation reveals
$$
|\Pnum (N_n (x_n^*) \geq 1) -[1-(1-\pi_{x_n^*})^n]| \leq O
(\frac{(\log n)^3}{n^{0.5+\vep_0}}).
$$

Summing the above results together and noting $\widetilde{\Pnum}
(N_n (x) \geq 1)=1-(1-\pi_x)^n$, we have proved the following lemma.
\begin{lem}\label{lem: A1}
Let $x_n^*=\left \lfloor \sqrt{\frac{Cn}{\log n}} \right \rfloor$ with $C=4+\frac{2}{\log 2}$.
Then
\begin{eqnarray*}
\sum_{x > x_n^*} |\Pnum (N_n (x) \geq 1) -[1-(1-\pi_x)^n]| &\leq& O
(n^{-0.2}),\\
|\Pnum (N_n (x_n^*) \geq 1) -[1-(1-\pi_{x_n^*})^n]| &\leq& O
(n^{-0.7}).
\end{eqnarray*}
\end{lem}

Now it is natural to see what happens for $|\Pnum (N_n (x) \geq 1)
-[1-(1-\pi_x)^n]|$ with $x \leq x_n^*$. The inequality $0 \leq \Pnum
(N_n (x) \geq 1) \leq 1$ is surely not sufficient for our purpose in
view of Lemma \ref{lem: main-tool}. Intuitively, we should have the
following result.
\begin{lem}\label{lem: A2}
$f (x):=\Pnum (N_n (x) \geq 1)$ is decreasing in $x$.
\end{lem}
But a mathematical rigorous proof is not so obvious. We postpone its
proof a bit later.

From Lemma \ref{lem: A2} we should have
$$
f (x_n^*) \leq f (x) \leq 1, \hbox{ for } x=1, 2, \cdots, x_n^*.
$$
And we already know from Lemma \ref{lem: A1} that
$$
f (x_n^*) \geq 1-(1-\pi_{x_n^*})^n-O (n^{-0.7}).
$$
Since $1/(C \log 2) =\frac{1}{2+4 \log 2} =\vep_0 =0.209529...$, we have for sufficiently
large $n$
$$
f (x_n^*) \geq 1-n^{-0.2},
$$
which implies the following
\begin{lem}\label{lem: A3}
Let $x_n^*$ be defined as above. Then for sufficiently large $n$
$$
|\Pnum (N_n (x) \geq 1) -[1-(1-\pi_{x})^n]| \leq n^{-0.2},
\qquad x=1, 2, \cdots, x_n^*.
$$
\end{lem}

In order to prove Lemma \ref{lem: A2}, we observe the following
important fact.
\begin{lem}\label{lem: important-obs}\textbf{(Comparison Lemma for Continued Fraction)}
Given two sequences of natural numbers $\tilde{x}=(x_k: 1 \leq k \leq n)$
and $\tilde{y}=(y_k: 1 \leq k \leq n)$. Suppose $x_k \geq y_k$ for $k=1,
\cdots, n$. Then $\mu (I (\tilde{x})) \leq \mu (I (\tilde{y}))$.
\end{lem}

\noindent \textit{Proof}. \;
In fact, let $\displaystyle \frac{p_k}{q_k}=[x_1, \cdots, x_k], \frac{\bar{p}_k}{\bar{q}_k}
=[y_1, \cdots, y_k]$ be irreducible fractions. Then
$$
\mu (I (\tilde{x})) = \int_0^1 \frac{\rmd \om}{(\log 2) \cdot (q_n +\om q_{n-1}) \cdot
[p_n+q_n +\om (p_{n-1}+q_{n-1})]}=: \int_0^1 \rho (\om; \tilde{x})\rmd \om
$$
and similar equation holds for $\mu (I (\tilde{y}))$. Then the condition in the lemma implies
$$
p_{n-1} \geq \bar{p}_{n-1}, \quad q_{n-1} \geq \bar{q}_{n-1}, \quad
p_n \geq \bar{p}_n, \quad q_n \geq \bar{q}_n.
$$
Therefore the densities satisfy $\rho (\om; \tilde{x}) \leq \rho
(\om; \tilde{y})$, which implies $\mu (I (\tilde{x})) \leq \mu (I
(\tilde{y}))$.
\qed

\begin{rem}
The proof of Lemma \ref{lem: important-obs} says more. Let $\tilde{x}=(x_1,
\cdots, x_m), \tilde{y}=(y_1, \cdots, y_n)$ be two natural number tuples
which may be of different length. Let
$$
\frac{p_m}{q_m} :=[x_1, \cdots, x_m], \quad
\frac{\bar{p}_n}{\bar{q}_n} :=[y_1, \cdots, y_m]
$$
be irreducible fractions. If
$$
\left[
    \begin{array}{cc}
p_{m-1}   &p_m\\
q_{m-1}\; &q_m
    \end{array}
\right]
\geq
\left[
    \begin{array}{cc}
\bar{p}_{n-1}   &\bar{p}_n\\
\bar{q}_{n-1}\; &\bar{q}_n
    \end{array}
\right],
$$
then $\mu (I (\tilde{x})) \leq \mu (I (\tilde{y}))$.
\end{rem}

\noindent \textit{Proof of Lemma \ref{lem: A2}}. \;
For any $x<y$, we would prove $f (x) \geq f (y)$. Noting that
$$
f (x):=\Pnum (N_n (x) \geq 1)=\Pnum (N_n (x) \geq 1, N_n
(y)=0)+\Pnum (N_n (x) \geq 1, N_n (y) \geq 1),
$$
we only need to prove for any $k=1, 2, \cdots, n$
$$
\Pnum (N_n (x)=k, N_n (y)=0) \geq \Pnum (N_n (y)=k, N_n (x)=0).
$$
But this is obvious in view of Lemma \ref{lem: important-obs} since
$$
\Pnum (N_n (x)=k, N_n (y)=0)=\sum_{***} \mu (I (x_1, \cdots, x_n)),
$$
where $\displaystyle \sum_{***}$ is over all tuples $(x_1, \cdots, x_n) \in \Nnum^n$
with $x_i \neq y, \forall i$ and $\displaystyle \sum_{i=1}^n 1_{\{x_i=x\}}=k$.
\qed

Summing the above estimations together, we have proved
\begin{lem}\label{lem: A4}
For sufficiently large $n$, $\displaystyle |\Enum R_n-\widetilde{\Enum} R_n| \leq O (n^{0.3})$.
And by Lemma \ref{lem: basic-facts},
$$
\Enum R_n=\sqrt{\frac{\pi n}{\log 2}} +\overline{O}
(n^{0.3}).
$$
\end{lem}
\subsection{Estimating $\mathrm{Var\,} (R_n)$}\label{sec:3.2}
Now we come to the estimation of $\mathrm{Var\,}
(R_n)$. We would need to estimate the probability
$$
\Pnum (N_n (x) \geq 1, N_n (y) \geq 1) \hbox{ for all } x \neq y.
$$

Similar to the above subsection, we need the following elemental
fact.
\begin{lem}\label{lem: Basic-2}
Let $\{A_k\}_1^n, \{B_k\}_1^n$ be two sequences of measurable sets
in a probability space $(\Om, \calF, \Pnum)$ which satisfies $A_k
\bigcap B_k=\emptyset, k=1, \cdots, n$. Then
\begin{equation}\label{eq: basic-eq-2}
\Pnum (N_A \geq 1, N_B \geq 1)=\sum_{k=2}^n (-1)^{k}
\sum_{\stackrel{r+s=k}{1 \leq r <k}} \sum_{\stackrel{i_1, \cdots, \,
i_r}{j_1, \cdots, \, j_s}} \Pnum (A_{i_1, \cdots, \, i_r}
B_{j_1, \cdots, \, j_s}),
\end{equation}
where in the above summation the indices $i_1, \cdots, i_r$ and
$j_1, \cdots, j_s$ are all distinct and both groups of indices are
in increasing order.
\end{lem}
One can prove eq. (\ref{eq: basic-eq-2}) directly using eq.
(\ref{eq: basic-eq'}).

Then we immediately derive an estimation
\begin{eqnarray}
\nonumber \Delta (x, y) &:=& |\Pnum (N_n (x) \geq 1, N_n (y) \geq
1)-\widetilde{\Pnum} (N_n (x) \geq 1, N_n (y) \geq 1)|\\
\label{eq: 3.7} &\leq& \sum_{k=2}^n \sum_{\stackrel{r+s=k}{1 \leq r <
k}} \sum_{\stackrel{i_1, \cdots, \, i_r}{j_1, \cdots, \, j_s}}
|\Pnum (A^x_{i_1, \cdots, \, i_r} A^y_{j_1, \cdots, \,
j_s})-\widetilde{\Pnum} (A^x_{i_1, \cdots, \, i_r} A^y_{j_1,
\cdots, \, j_s})|.
\end{eqnarray}

We first make the estimation for $x>y \geq x_n^*$.
Similar calculation yields:
\begin{lem}\label{lem: B1}
For sufficiently large $n$, $\displaystyle \sum_{x>y \geq x_n^*} \Delta (x, y) \leq O (\sqrt{n})$.
\end{lem}

For $x<y \leq x_n^*$, we have
\begin{eqnarray}
\nonumber && \Pnum (N_n (x) \geq 1, N_n (y) \geq 1)-\widetilde{\Pnum} (N_n (x)
\geq 1, N_n (y) \geq 1) \\
\nonumber &\leq& |\Pnum (N_n (y) \geq 1)-\widetilde{\Pnum} (N_n (y) \geq
1)|+\widetilde{\Pnum} (N_n (y) \geq 1) -\widetilde{\Pnum} (N_n (x)
\geq 1, N_n (y) \geq 1)\\
\label{eq:estimate4var} &=&|\Pnum (N_n (y) \geq 1)-\widetilde{\Pnum} (N_n (y) \geq
1)|+\widetilde{\Pnum} (N_n (x) =0, N_n (y) \geq 1)\\
\nonumber &\leq& O (n^{-0.2})+\widetilde{\Pnum} (N_n (x) =0)= O
(n^{-0.2})+(1-\pi_x)^n \\
\nonumber &\leq& O (n^{-0.2})+(1-\pi_{x^*_n})^n=O (n^{-0.2}).
\end{eqnarray}
Hence we have the following lemma.
\begin{lem}\label{lem: B1'}
For sufficiently large $n$,
$$
\sum_{x<y\leq x_n^*} [\Pnum (N_n (x) \geq 1, N_n (y) \geq
1)-\widetilde{\Pnum} (N_n (x) \geq 1, N_n (y) \geq 1)] \leq O
(n^{0.8}).
$$
\end{lem}

Similarly, for $x<x_n^* \leq y$, noting eq. (\ref{eq:estimate4var}) we have
\begin{eqnarray*}
&& \Pnum (N_n (x) \geq 1, N_n (y) \geq
1)-\widetilde{\Pnum} (N_n (x) \geq 1, N_n (y) \geq 1)\\
&\leq& |\Pnum (N_n (y) \geq 1)-\widetilde{\Pnum} (N_n (y) \geq 1)| +\widetilde{\Pnum} (N_n
(x) =0, N_n (y) \geq 1)\\
&\leq& |\Pnum (N_n (y) \geq 1)-\widetilde{\Pnum} (N_n (y) \geq 1)| +[(1-\pi_x)^n-(1-\pi_x-\pi_y)^n]
\end{eqnarray*}
which implies
\begin{eqnarray*}
&& \sum_{x<x_n^* \leq y} [\Pnum (N_n (x) \geq 1, N_n (y) \geq
1)-\widetilde{\Pnum} (N_n (x) \geq 1, N_n (y) \geq 1)]\\
&\leq& O (\sqrt{\frac{n}{\log n}}) \cdot O
(n^{-0.2})+\sum_{x<x_n^* \leq y} [(1-\pi_x)^n
-(1-\pi_x-\pi_y)^n]\\
&\leq& O (n^{0.3})+\sum_{x<x_n^*} (1-\pi_x)^n \cdot \sum_{y
\geq x^*_n} [1 -(1-\frac{\pi_y}{1-\pi_x})^n]\\
&\leq& O (n^{0.3})+\sum_{x<x_n^*} (1-\pi_x)^n \cdot \sum_{y \geq x^*_n} [1 -(1-2
\cdot \pi_y)^n] \hbox{ (noting } \pi_x \leq \pi_{_1}=0.4150...\hbox{)}\\
&=& O (n^{0.3})+ \sum_{x<x_n^*} (1-\pi_x)^n \cdot O (\sqrt{n})  \qquad \qquad (*)\\
&\leq&O (n^{0.3})+ O (\sqrt{n}) \cdot O (\sqrt{\frac{n}{\log n}}) \cdot O
(n^{-0.2095...}) \leq O (n^{0.8}).
\end{eqnarray*}
We note here that, in the step $(*)$ we have exploited the estimating technique
developed in \cite{CXY}. Thus we have the following lemma.
\begin{lem}\label{lem: B1'}
For sufficiently large $n$,
$$
\sum_{x<x_n^* \leq y} [\Pnum (N_n (x) \geq 1, N_n (y) \geq
1)-\widetilde{\Pnum} (N_n (x) \geq 1, N_n (y) \geq 1)] \leq
O (n^{0.8}).
$$
\end{lem}

Note that, by Lemma \ref{lem: basic-facts},
$$
\mathrm{Var}_{\widetilde{\Pnum}} (R_n) \leq \widetilde{\Enum}
R_n =O (\sqrt{n})
$$
and in view of the above estimations
\begin{eqnarray*}
\mathrm{Var\,} (R_n)-\mathrm{Var}_{\widetilde{\Pnum}} (R_n)
&=& \Enum R_n^2-\widetilde{\Enum} R_n^2+(\widetilde{\Enum} R_n)^2
-(\Enum R_n)^2\\
&=& \Bigl[ \Enum R_n-\widetilde{\Enum} R_n \Bigr] +\Bigl[
(\widetilde{\Enum} R_n)^2 -(\Enum R_n)^2 \Bigr] \\
&&+2 \sum_{x<y} \Bigl[ \Pnum (N_n (x) \geq 1, N_n (y) \geq 1)- \widetilde{\Pnum} (N_n (x)
\geq 1, N_n (y) \geq 1) \Bigr]\\
&\leq& O (n^{0.3})+O (n^{0.3} \cdot \sqrt{n}) +O
(n^{0.8})=O (n^{0.8}),
\end{eqnarray*}
we have the following lemma.
\begin{lem}\label{lem: B2}
For sufficiently large $n$,
$$
\mathrm{Var\,} (R_n) \leq O (n^{0.8}) \leq O
(1) \cdot (\Enum R_n)^{2-\delta} \quad \hbox{ with } \; \delta=0.4.
$$
\end{lem}

Now by Lemma \ref{lem: main-tool}, we have for $\mu$-a.e. (and hence
for Lebesgue almost all) $\om \in (0,1)$
$\displaystyle \lim_{n \to +\infty} \frac{R_n (\om)}{\Enum R_n}=1$, i.e.,
$$
\lim_{n \to +\infty} \frac{R_n
(\om)}{\sqrt{n}}=\sqrt{\frac{\pi}{\log 2}}.
$$

\subsection{More estimations for $R_{n, \, k+}$ with $k \geq 2$}\label{sec:3.3}
The estimation of $\Enum R_{n, \, k+}$ and $\mathrm{Var\,} (R_{n, \,
k+})$ follows almost the same line as in the previous subsections
except that we need some additional treatments.

First we would need new equations in the places of eq. (\ref{eq:
basic-eq}) and (\ref{eq: basic-eq-2}), which we state as the
following lemmas without proof.

\begin{lem}\label{lem: Basic-3}
Let $\{A_i\}_1^n$ be a sequence of measurable sets in a probability
space $(\Om, \calF, \Pnum)$. Then for any $1 \leq k \leq n$
$$
\Pnum (N_A \geq k)=\sum_{r=k}^n (-1)^{r-k} \cdot C_{r-1}^{k-1} \cdot
\sum_{1 \leq i_1 <\cdots< i_r \leq n} \Pnum (A_{i_1, \cdots, \,
i_r}).
$$
\end{lem}

\begin{lem}\label{lem: Basic-4}
Let $\{A_i\}_1^n, \{B_i\}_1^n$ be two sequences of measurable sets
in a probability space $(\Om, \calF, \Pnum)$ such that $A_i \cap B_i
=\emptyset$ for all $i$. Then for any $1 \leq k \leq n$
$$
\Pnum (N_A \geq k, N_B \geq k)=\sum_{r=2 k}^n (-1)^{r} \cdot
\sum_{\stackrel{a+b=r}{a, \, b \geq k}} C_{a-1}^{k-1} \cdot
C_{b-1}^{k-1} \cdot \sum_{\stackrel{i_1, \cdots, \, i_a}{j_1,
\cdots, \, j_b}} \Pnum (A_{i_1, \cdots, \, i_a} \bigcap B_{j_1,
\cdots, \, j_b}),
$$
where in the above summation the two group of indices $i_1, \cdots,
\, i_a$ and $j_1, \cdots, \, j_b$ are all distinct and both in
increasing order.
\end{lem}

We believe that there is monotonicity in the function
\begin{equation}
f_k (x) :=\Pnum (N_n (x) \geq k)
\end{equation}
for all (fixed) $n \geq k \geq 2$ as Lemma \ref{lem: A2} states. But
we cannot give a rigorous and relatively easy proof for such result. Therefore we
present an alternative treatment here.
\begin{lem}\label{lem: C-3}
For all $1 \leq x \leq x^*_n$ and sufficiently large $n$,
$\displaystyle \Pnum (N_n (x) \geq k) \geq 1+O (n^{-0.2})$. Hence
$$
|\Pnum (N_n (x) \geq k)- \widetilde{\Pnum} (N_n (x) \geq k)|
\leq O (n^{-0.2}).
$$
\end{lem}
\noindent \textit{Proof}. \;
Assume $k \geq 2$. Let $L_n$ be defined as above and let $s_1 :=0$. Put
\begin{eqnarray*}
A_n &:=& \left \lfloor \frac{n- (k-1) \cdot L_n}{k} \right \rfloor, \\
t_j &:=& (j-1) \cdot L_n+j \cdot A_n, \quad s_{j+1}=j \cdot (L_n+
A_n), j=1, \cdots, k.
\end{eqnarray*}
Clearly $t_k \leq n$ and $0=s_1<t_1<s_2<t_2<\cdots<s_k<t_k \leq n$. We will write
$$
N_x (\Delta):=\sum_{i \in \Delta} 1_{\{a_i (\om)=x\}}
$$
for any interval $\Delta \subset \Nnum$, which is the visiting number at $x$ with
times $n \in \Delta$ such that $a_n (\om)=x$. Then obviously
$$
\Pnum (N_n (x) \geq k) \geq \Pnum (N_x ((s_j, t_j]) \geq 1 \hbox{
for } j=1, \cdots, k).
$$
Hence
\begin{eqnarray*}
\Pnum (N_n (x) \geq k) &\geq& (1-O (q^{L_n}))^{k-1} \cdot \prod_{i=1}^k
\Pnum (N_x ((s_j, t_j]) \geq 1)\\
&=& (1-O (q^{L_n}))^{k-1} \cdot \Bigl[ \Pnum (N_x ((0, A_n]) \geq 1) \Bigr]^k\\
&\geq& [1-O (n^{-4})]^{k-1} \cdot [1-O (n^{-0.2})]^k\\
&=& 1+O (n^{-0.2}).
\end{eqnarray*}
We have the same bound for $\widetilde{\Pnum} (N_n (x) \geq k)$.
\qed

With the help of the above lemmas one can prove
$$
\lim_{n \to +\infty} \frac{R_{n, \, k+} (\om)}{\Enum R_{n, \, k+}}=1
$$
following the same way indicated above; the details are omitted
here. By Lemma \ref{lem: basic-facts}
$$
\lim_{n \to +\infty} \frac{\Enum R_{n, \, k+}}{\Enum R_n}=r_{_{k+}},
\quad \forall k \geq 1.
$$
Therefore our main theorem follows.

\section{Proof of Theorem \ref{thm: main2}}\label{sec:4}

\subsection{Some elementary facts on continued fractions}
In this subsection, we collect some elementary properties shared by
continued fractions that will be used later.

For any $ n\geq 1$ and $(a_1, a_2, \cdots, a_n) \in \mathbb{N}^n$,
let $q_n(a_1, a_2, \cdots, a_n)=q_n$ be defined by (\ref{f1}). Then
we have
\begin{prop}\label{p2} (See \cite{Wu}.)
For any $ n\geq 1$ and $1 \leq k \leq n$,
\begin{eqnarray}\label{f4}
\frac{a_k+1}{2} \leq \frac{q_n(a_1, a_2, \cdots, a_n)}{q_{n-1}(a_1,
\cdots, a_{k-1}, a_{k+1},\cdots, a_n)} \leq a_k+1.
\end{eqnarray}
\end{prop}

For any positive integer $B\geq 2$, Let $E_B$ be the set of
continued fractions each of whose partial quotients is between $1$
and $B$, i.e.
$$
E_B=\Big\{ \om \in [0,1):\ 1 \leq a_n(\om) \leq B, \ {\rm{for\
all}}\ n\geq 1\ \Big\}.
$$
I. J. Good \cite{Go} proved the following result.
\begin{prop}\label{p3} (See \cite{Go}.)
For any $n \geq 1$, let $\sigma_n$ be the unique root of
$$
\sum\limits_{1\leq a_1, \, a_2, \, \cdots, \, a_n \leq
B}\frac{1}{q_n(a_1, a_2, \cdots, a_n)^{2s}}=1.
$$
Then
$$
\dim_H E_B=\lim\limits_{n \to \infty} \sigma_n.
$$
Moreover, $\lim\limits_{B \to \infty} \dim_H E_B=1$.
\end{prop}
\subsection{Non-autonomous conformal iterated function systems}
In this part, we present the construction and some basic properties
of non-autonomous conformal \textbf{iterated function system} (abbr. \textbf{IFS}) which was
introduced quite recently by L. Rempe-Gillen and M. Urbanski in
\cite{RU}. It is a variant of the construction studied in
\cite{HRWW} and \cite{Wen}.

Fix a compact set $X \subset \Rnum^d$ with
$\overline{int(X)}=X$ such that $\partial X$ is smooth or $X$ is
convex. Given a conformal map $\varphi:\ X \to X$, we denote by
$D\varphi(x)$ the derivative of $\varphi$ at $x$, and we denote by
$|D\varphi(x)|$ the operator norm of the differential. Put
$$
\|D\varphi\|=\sup\{|D\varphi(x)|:\ x \in X\},\
|||D\varphi|||=\inf\{|D\varphi(x)|:\ x \in X\}.
$$

For any $n\geq 1$, let $ I^{(n)}$ be a (finite or countable
infinite) index set. For each $i\in I^{(n)}$, there is a conformal
map $\varphi^{(n)}_i:\ X \to X$, and write
$\Phi^{(n)}=\Big\{\varphi^{(n)}_i:\ i\in I^{(n)}\Big\}$.
\begin{defn}\label{d1}
We call $\Phi=\Big\{\Phi^{(1)},\ \Phi^{(2)},\
\Phi^{(3)}\cdots\Big\}$ a non-autonomous conformal IFS on the set $X$ if the following holds.
\begin{itemize}
  \item[{\rm (a)}] Open set condition: We have
$$
\varphi_i^{(n)}(int(X))\bigcap \varphi_j^{(n)}(int (X))=\emptyset.
$$
for all $n \in \Nnum$ and all distinct indices $i, j \in
I^{(n)}$.
  \item[{\rm (b)}] Conformality: There exists an open connected set $V\supset X$
such that for each $n\geq 1$ and $i \in I^{(n)}$, $\varphi_i^{(n)}$
can be extended to a $C^1$ conformal diffeomorphism of $V$ into $V$.
  \item[{\rm (c)}] Bounded distortion: There exists a constant $K \geq 1$ such that,
for any $k \leq l$ and any $i_k,\ i_{k+1}, \cdots, i_l$ with
$i_j \in I^{(j)}$, the map
$\varphi=\varphi_{i_k}^{(k)}\circ\cdots\circ\varphi_{i_l}^{(l)}$
satisfies
$$
|D\varphi(x)|\leq K |D\varphi(y)|
$$
for all $x,\ y \in V$.
  \item[{\rm (d)}] Uniform contraction: There exists a constant $\eta<1$ such that
$$
|D \varphi(x)|\leq \eta^m
$$
for all sufficiently large $m$, all $x \in X$ and all
$\varphi=\varphi_{i_j}^{(j)}\circ\cdots\circ\varphi_{i_{j+m}}^{(j+m)}$,
where $j \geq 1$ and $i_k \in I^{(k)}$.
\end{itemize}
\end{defn}

For any $0<m \leq n <\infty$, write
$$
I^n :=\prod_{j=1}^n I^{(j)},\ I^\infty :=\prod_{j=1}^{\infty} I^{(j)},\
I^{m, n} :=\prod_{j=m}^n I^{(j)},\ {\rm{and}}\ I^{m,
\infty} :=\prod_{j=m}^{\infty} I^{(j)}.
$$
If $\tilde{i}=i_m i_{m+1}\cdots i_n \in I^{m, n}$, write $\varphi^{m,
n}_{\tilde{i}}=\varphi_{i_1}^{(m)}\circ\cdots\circ\varphi_{i_n}^{(n)}$.
When $m=1$, we also abbreviate $\varphi_{\tilde{i}} :=\varphi^n_{\tilde{i}}
:=\varphi^{1, n}_{\tilde{i}}$.

For any $n \geq 1$ and $i \in I^n$, let $X_i=\varphi_i (X)$.
The limit set (or attractor) of $\Phi$ is defined as
\begin{equation}
J :=J(\Phi):=\bigcap_{n=1}^{\infty}\bigcup_{i \in I^n} X_i.
\end{equation}
\begin{defn}\label{d2}
For any $t \geq 0$ and $n \in \Nnum$, we define
$$
Z_n(t)=\sum\limits_{i \in I^n}\|D\varphi_i\|^t,
$$
the upper and lower pressure functions are defined as follows:
$$
\underline{P}(t)=\varliminf\limits_{n \to \infty}\frac{1}{n} \log
Z_n(t),\ \ \overline{P}(t)=\varlimsup\limits_{n \to \infty}\frac{1}{n}
\log Z_n(t).
$$
\end{defn}

L. Rempe-Gillen and M. Urbanski \cite{RU} proved the following
results.
\begin{prop}\label{p4}
If $\lim\limits_{n \to \infty} \frac{1}{n} \log \sharp I^{(n)}=0$,
then
\begin{eqnarray*}
\dim_H J&=&\sup\{t \geq 0:\ \underline{P}(t)>0\}=\inf\{t \geq 0:\ \underline{P}(t)<0\}\\
&=&\sup\{t\geq 0:\ Z_n(t) \to \infty\}.
\end{eqnarray*}
\end{prop}
\begin{prop}\label{p5}
Suppose that both limits
$$
a=:\lim\limits_{n \to \infty} \frac{1}{n} \log \sharp I^{(n)}
$$
and
$$
b=:\lim\limits_{n \to \infty, j\in I^{(n)}} \frac{1}{n} \log
\Big(\frac{1}{\|D\varphi_j^{(n)}\|}\Big)
$$
exist and positive finite. Then $\dim_H J=\frac{a}{b}$.
\end{prop}
\subsection{Proof of Eq. (\ref{eq: dim-1}) and (\ref{eq: dim-2})}
{\sc Proof of Eq. (\ref{eq: dim-1})}:
We only prove for $c=1$, $E(\beta):= E(\beta, c)$ is of Hausdorff
dimension $1$, the others can be proved in the same way.

For any $\vep>0$, by Proposition \ref{p3}, we can find a positive
integer $B$ with $\dim E_B> 1-\frac{\vep}{4}$, and an integer $N_0>0$ such that
$\sigma_n>1-\frac{\vep}{2}$ for any $n \geq N_0$ .

Put $\gamma :=\frac{1}{\beta}>1$. Define a subset $E_B(\beta)$ of $E(\beta)$ as follows.
\begin{eqnarray*}
E_B(\beta)&=&\Big\{\om \in [0,1):\ a_{_{\left \lfloor k^{\gamma} \right \rfloor}} (\om)=k\
{\rm{for\ all}}\  k\geq 1, \ {\rm{and}}\ \\
& & \ 1\leq a_n (\om) \leq B, \ {\rm{for\ all \ other}}\ n\geq 1\
\Big\}.
\end{eqnarray*}
It is direct to check that
$$
E_B(\beta)\subset E(\beta).
$$

Define $\Phi=\{\Phi^{(1)},\ \Phi^{(2)}, \cdots\}$ as follows. If
$n=\left \lfloor k^{\gamma} \right \rfloor$ for some $k\geq 1$ (noting that such $k$ is unique!), let
$\Phi^{(n)}=\{\frac{1}{x+k}\}$; let $\Phi^{(n)}=\Big\{\frac{1}{x+1},\ \frac{1}{x+2}, \cdots,
\frac{1}{x+B}\Big\}$ for other $n$. Then $\Phi$ is a non-autonomous conformal IFS, and $E_B(\beta)$ is the associated
limit set. $I^{(n)}=\{k\}$, if $n=\left \lfloor k^{\gamma} \right \rfloor$ for some $k \geq 1$,
and $I^{(n)}=\{1,\ 2,\cdots, B\}$ for others. Thus $\sharp I^{(n)}=1$, if
$n=\left \lfloor k^{\gamma} \right \rfloor$ for some $k \geq 1$, and
$\sharp I^{(n)}=B$ for others. By Proposition \ref{p4}, we have
$$
\dim E_B(\beta)=\sup\{t\geq 0:\ Z_n(t) \to \infty\}.
$$

Now for any $n \geq 1$, choose $k \geq 1$ such that
$\left \lfloor k^{\gamma} \right \rfloor \leq n <\left \lfloor (k+1)^{\gamma} \right \rfloor$. By
(\ref{f4}), we have
\begin{eqnarray*}
&& Z_n(1-\vep)=\sum\limits_{\tilde{i} \in I^n}\|D\varphi_{\tilde{i}} \|^{(1-\vep)}\\
&\geq& \prod_{j=1}^{k-1} \Big( \sum \limits_{\tilde{i} \in I^{\left \lfloor j^{\gamma} \right \rfloor,
\left \lfloor (j+1)^{\gamma} \right \rfloor-1}} |||D\varphi_{\tilde{i}}^{\left \lfloor j^{\gamma} \right \rfloor,
\left \lfloor (j+1)^{\gamma} \right \rfloor-1} |||^{1-\vep} \Big) \cdot \sum\limits_{\tilde{i}
\in I^{\left \lfloor k^{\gamma} \right \rfloor, n}} |||D\varphi_{\tilde{i}}^{\left \lfloor k^{\gamma}
\right \rfloor, n} |||^{1-\vep}.
\end{eqnarray*}
For any $j \geq 1$, write $l(j)=\left \lfloor (j+1)^{\gamma} \right \rfloor-\left \lfloor j^{\gamma}
\right \rfloor$. Notice that when $j$ is large enough such that $l(j)-1\geq N_0$, we have
\begin{eqnarray}\label{f5}
&&\sum\limits_{\tilde{i} \in I^{\left \lfloor j^{\gamma} \right \rfloor, \left \lfloor (j+1)^{\gamma} \right \rfloor -1}}
|||D\varphi_{\tilde{i}}^{\left \lfloor j^{\gamma} \right \rfloor, \left \lfloor (j+1)^{\gamma} \right \rfloor -1}
|||^{1-\vep}\nonumber\\
&=&\sum\limits_{1\leq a_1, a_2,\cdots, a_{l(j)-1} \leq B} \Big(
\frac{1}{q_{l(j)} (j, a_1, a_2, \cdots, a_{l(j)-1})+q_{l(j)-1}(j,
a_1, a_2, \cdots, a_{l(j)-2})} \Big) ^{2(1-\vep)}\nonumber\\
&\geq& \frac{1}{2^{2(1-\vep)}} \sum\limits_{1\leq a_1, a_2,\cdots,
a_{l(j)-1} \leq B}
\Big(\frac{1}{q_{l(j)}(j, a_1, a_2, \cdots, a_{l(j)-1})}\Big)^{2(1-\vep)}\nonumber\\
&\geq&\frac{1}{(2(j+1))^{2(1-\vep)}} \sum\limits_{1\leq a_1,
a_2,\cdots, a_{l(j)-1}
\leq B}\Big(\frac{1}{q_{l(j)-1}(a_1, a_2, \cdots, a_{l(j)-1})} \Big)^{2(1-\vep)}\nonumber\\
&\geq& \frac{1}{(2(j+1))^{2(1-\vep)}} \sum\limits_{1\leq a_1,
a_2,\cdots, a_{l(j)-1} \leq B} \Big(\frac{1}{q_{l(j)-1}(a_1, a_2,
\cdots, a_{l(j)-1})}\Big)
^{2\sigma_{l(j)-1}-\vep}\nonumber\\
&\geq &  \frac{1}{(2(j+1))^{2(1-\vep)}} \cdot
2^{\frac{l(j)-1}{2}\vep} \geq \frac{1}{(2(j+1))^{2}} \cdot
2^{\frac{(j+1)^{\gamma}-j^{\gamma}-3}{2}\vep}.
\end{eqnarray}

\noindent In the similar way, we have, if $n-\left \lfloor k^{\gamma} \right \rfloor-1
\geq N_0$,
\begin{equation}\label{f6}
\sum\limits_{\tilde{i} \in I^{\left \lfloor k^{\gamma} \right \rfloor, n}}
|||D\varphi_{\tilde{i}}^{\left \lfloor k^{\gamma} \right \rfloor, n} |||^{1-\vep} \geq
\frac{1}{(2(k+1))^{2}} \cdot
2^{\frac{n-k^{\gamma}-2}{2}\vep}.
\end{equation}

\noindent If $n-[k^{\gamma}]-1 < N_0$,
\begin{equation}\label{f7}
\sum\limits_{\tilde{i} \in I^{[k^{\gamma}], n}}
|||D\varphi_{\tilde{i}}^{[k^{\gamma}], n} |||^{1-\vep} \geq
\frac{1}{(2(k+1))^{2}} \cdot \frac{1}{(B+1)^{N_0}}.
\end{equation}
Combining (\ref{f5}), (\ref{f6}) and (\ref{f7}), we have
$\lim\limits_{n \to \infty} Z_n(1-\vep)=\infty$. This implies
$\dim E_B(\beta)\geq 1-\vep$. Since $\vep$ is arbitrary, we finish
the proof of Eq. (\ref{eq: dim-1}).
\qed
\bigskip

\noindent {\sc Proof of Eq. (\ref{eq: dim-2})}:
We divide the proof into two parts.

\bigskip

{\sc Upper Bound}: For any $n\geq 1$ and $(a_1,\cdots, a_n)\in
\Nnum^n$, let $R_n(a_1, a_2, \cdots, a_n)$ be the number of
distinct ones among $a_1, a_2, \cdots, a_n$.

For any $0<\vep<c$, let $t=\frac{1}{2}+\vep$ and
$s=\frac{1}{2}+\frac{\vep}{2}$. For any $n\geq 1$, let
$$
\Lam_n=\Big\{(a_1,\cdots, a_n)\in \Nnum^n:\ R_n(a_1, a_2,
\cdots, a_n)\geq (c-\vep)n \Big\}.
$$
Then
$$
F (c) \subset \bigcup_{N=1}^{\infty} \bigcap_{n=N}^{\infty}
\bigcup_{(a_1,\cdots, a_n) \in \Lam_n} I(a_1,\cdots, a_n).
$$

\noindent For any $N\geq 1$,
\begin{eqnarray*}
&&{H}^{t}(\bigcap_{n=N}^{\infty}\bigcup_{(a_1,\cdots, a_n) \in
\Lam_n} I(a_1,\cdots, a_n))\leq \varliminf\limits_{n \to \infty}
\sum\limits_{(a_1,\cdots, a_n) \in \Lam_n}|I(a_1,\cdots, a_n)|^{t}\\
&\leq& \varliminf\limits_{n \to \infty} \sum\limits_{(a_1,\cdots,
a_n) \in \Lam_n} (a_1\cdot a_2 \cdots a_n)^{-2t}\\
&\leq&\varliminf\limits_{n \to \infty}
\frac{1}{\Big(\left \lfloor (c-\vep) n \right \rfloor!\Big)
^{\vep}}\sum\limits_{(a_1,\cdots, a_n) \in \Lam_n} (a_1\cdot a_2 \cdots a_n)^{-2s}\\
&\leq&\varliminf\limits_{n \to \infty}
\frac{1}{\Big( \left \lfloor (c-\vep) n \right \rfloor!\Big)
^{\vep}}\sum\limits_{(a_1,\cdots, a_n) \in \Nnum^n} (a_1\cdot a_2 \cdots a_n)^{-2s}\\
&=&\varliminf\limits_{n \to \infty} \frac{1}{\Big(\left \lfloor (c-\vep) n \right \rfloor!\Big)^{\vep}} \cdot
[\zeta (1+\vep)]^n  \hbox{ (where } \zeta (\cdot) \hbox{ is Riemann's zeta function)}\\
&=& 0.
\end{eqnarray*}
We finish the proof of {\sc upper bound}.

\bigskip

{\sc Lower Bound}: For any $c \in (0,1]$ given, let
\begin{eqnarray*}
G &:=& \Big\{ \om \in [0,1): 2^n \leq a_k (\om) <2^{n+1} \hbox{ if }
k=\left \lfloor n/c \right \rfloor \hbox{ for some } n \geq 1, \\
&& \hbox{ and for other } k, a_k (\om)=1 \Big\}.
\end{eqnarray*}
It is direct to check that
$$
G\subset F (c).
$$

Define $\Phi=\{\Phi^{(1)},\ \Phi^{(2)}, \cdots\}$ as follows. For
any $n \geq 1$, let $K_n=\left \lfloor n/c \right \rfloor, K^\prime_n =\left \lfloor (n+1)/c \right \rfloor$ and put
$$
\Phi^{(n)}=\Big\{[1, \cdots, 1, k+x]: \ 2^n \leq k<2^{n+1}\Big\},
$$
where the $1$'s in the continued fractional function formula $[1,
\cdots, 1, k+x]$ appear exactly $K^\prime_n-K_n-1$ (possibly
$K^\prime_n-K_n-1=0$) times. Then $\Phi$ is a non-autonomous conformal
IFS, and $G$ is the associated limit set. It is
easy to check that
$$
\lim\limits_{n \to \infty} \frac{1}{n} \log \sharp I^{(n)}=\log 2
$$
and
$$
\lim\limits_{n \to \infty, j\in I^{(n)}} \frac{1}{n} \log
\Big(\frac{1}{\|D\varphi_j^{(n)}\|}\Big)=2\log 2.
$$

By Proposition \ref{p5}, we have
$$
\dim_H G=\frac{1}{2},
$$
and this finish the proof of {\sc lower bound}.
\qed

\noindent{\sl \textbf{Acknowledgements} \quad} {The second author would like to thank
Prof. De-Jun Feng for telling him in early 2011 the strong mixing property for continued
fractions mentioned in the introduction part; He also thanks Mr. Peng Liu for the simulation
mentioned there. He is indebt to F. Ledrappier and Miss Jie Shen for helpful discussions.
This work is partially supported by NSFC (No. 11225101, No. 11171124, No. 10701026 and No.
11271077) and the Laboratory of Mathematics for Nonlinear Science, Fudan University.}





\bibliographystyle{elsarticle-num}

\begin{thebibliography}{00}


%
\bibitem{Billingsley} Billingsley, P.:
{\it Ergodic Theory and Information}, John Wiley, New York, 1965. MR0192027

\bibitem{CXY} Chen, X.-X.; Xie, J.-S.; Ying, J.-G.:
{\em Range-Renewal Processes: SLLN, Power Law and Beyonds},
preprint. arXiv:1305.1829

\bibitem{Derriennic} Derriennic, Y.:
{\em Quelques applications du th\'{e}or\`{e}me ergodique
sous-additif}. (French. English summary) Conference on Random Walks
(Kleebach, 1979) (French), pp. 183--201, 4, Ast\'{e}risque,
\textbf{74}, Soc. Math. France, Paris, 1980. MR0588163


\bibitem{DE50} Dvoretzky, A.; Erd\"{o}s, P.:
{\em Some problems on random walk in space}. Proceedings of the
Second Berkeley Symposium on Mathematical Statistics and
Probability, 1950. pp. 353--367. University of California Press,
Berkeley and Los Angeles, 1951. MR0047272


\bibitem{ET} Erd\"{o}s, P.; Taylor, S. J.:
{\em Some problems concerning the structure of random walk paths},
Acta Math. Acad. Sci. Hungar. \textbf{11} (1960), 137--162. MR0121870

\bibitem{Fa} Falconer, K. J.:
{\it Fractal Geometry, Mathematical Foundations and Application},
Wiley, 1990. MR1102677

\bibitem{Go} Good, I. J.:
{\em The fractional dimensional theory of continued fractions},
Proc. Camb. Philos. Soc., \textbf(37) (1941), 199--228. MR0004878

\bibitem{HRWW} Hua, S.; Rao, H.; Wen, Z.-Y.; Wu, J.:
{\em On the structures and dimensions of Moran sets}, Sci. China
Ser. A 43 (2000), no. 8, 836-852. MR1799919

\bibitem{Hi} Hirst, K. E.:
{\em A problem in the fractional dimension theory of continued fractions},
Quart. J. Math. Oxford Ser. \textbf(21) (1970), 29-35. MR0258778

\bibitem{I-K} Iosifescu, M.; Kraaikamp, C.:
{\it  Metrical Theory of Continued Fractions}, Mathematics and its
Applications, 547. Kluwer Academic Publishers, Dordrecht, 2002. MR1960327

\bibitem{Ja28} Jarn\'{i}k, V.:
{\em Zur metrischen Theorie der diophantine Approximationen},
Prace Mat.-Fiz., \textbf{36} (1928), 91--106.

\bibitem{Kh63} Khintchine, A. Ya.:
{\it Continued Fractions},
(Translated by Peter Wynn.) P. Noordhoff, Ltd., Groningen 1963 iii+101 pp. MR0161834

\bibitem{Kh64} Khintchine, A. Ya.:
{\it Continued Fractions},
The University of Chicago Press, Chicago, Ill.-London 1964 xi+95 pp. MR0161833

\bibitem{Kingman68} Kingman, J. F. C.:
{\em The ergodic theory of subadditive stochastic processes}, J.
Roy. Statist. Soc. Ser. B \textbf{30} (1968), 499--510. MR0254907

\bibitem{Kingman73} Kingman, J. F. C.:
{\em Subadditive ergodic theory}, Ann. Probability \textbf{1}
(1973), 883--909. MR0356192

\bibitem{Kingman76} Kingman, J. F. C.:
{\em Subadditive processes}, \'{E}cole d'\'{E}t\'{e} de
Probabilit\'{e}s de Saint-Flour, V-1975, pp. 167--223. Lecture Notes
in Math., Vol. \textbf{539}, Springer, Berlin, 1976. MR0438477

\bibitem{Kuzmin} Kuzmin, R. O.:
{\em On a problem of Gauss},
Dokl. Akad. Nauk SSSR Ser. A, 375--380. [Russian; French version in {\it
Atti Congr. Internaz. Mat.} (Bologna, 1928), Tomo VI, 83--89. Zanichelli,
Bologna, 1932]

\bibitem{LWWX} Li, B; Wang, B.-W.; Wu, J.; Xu, J.:
{\em The shrinking target problem in the dynamical system of continued fractions},
Proc. London Math. Soc., to appear.

\bibitem{Levy29} L\'{e}vy, P.:
{\em Sur les lois de probabilit\'{e} dont d\'{e}pendent les quotient complets
et incomplets d'une fraction continue},
Bull. Soc. Math. France \textbf{57} (1929), 178--194. MR1504948

\bibitem{Levy54} L\'{e}vy, P.:
{\em Th\'{e}orie de l'addition des variables al\'{e}atoires},
2\`{e}me \'{e}dition. Gauthier-Villars, Paris. (1\`{e}re \'{e}dition 1937). MR1509369

\bibitem{Lu} T. L\'uczak, {\em On the fractional dimension of sets
of continued fractions}, Mathematika 44 \textbf(1) (1997), 50-53. MR1464375

\bibitem{MU1} Mauldin, R. D.;Urba\'nski, M.:
{\em Dimensions and measures in infinite iterated function systems},
Proc. London Math. Soc. \textbf(73) (1996) 105-154.  MR1387085

\bibitem{MU2} Mauldin, R. D.; Urba\'nski, M.:
{\em Conformal iterated function systems with applications to the geometry of continued
fractions}, Trans. Amer. Math. Soc. 351 \textbf(12) (1999), 4995-5025. MR1487636

\bibitem{PW} Pollicott, M.; Weiss, H.:
{\em Multifractal analysis of Lyapunov exponent for continued fraction and
Manneville-Pomeau transformations and applications to Diophantine approximation},
Comm. Math. Phys. \textbf(207) (1999) 145-171. MR1724859

\bibitem{RU} Rempe-Gillen, L.; Urba\'nski, M.:
{\em Non-autonomous conformal iterated function systems and
Moran-set constructions}, arXiv:1210.7469.

\bibitem{Revesz} Revesz, P.:
{\it Random Walk in Random and Non-Random Environments},
Second edition. World Scientific Publishing Co. Pte. Ltd., Hackensack, NJ, 2005.
xvi+380 pp. ISBN: 981-256-361-X. MR2168855

\bibitem{Szusz} Sz\H{u}sz, P.:
{\em \"{U}ber einen Kusminschen Statz},
Acta Math. Acad. Sci. Hungar. \textbf{12} (1961), 447--453. MR0150121

\bibitem{WW} Wang, B.-W.; Wu, J.:
{\em Hausdorff dimension of certain sets arising in continued fraction expansions},
Adv. Math. \textbf(218) (2008) 1319-1339. MR2419924

\bibitem{Wen} Wen, Z.-Y.:
{\em Moran sets and Moran classes},
Chinese Sci. Bull. 46 (2001), no. 22, 1849--1856. MR1877244

\bibitem{Wirsing} Wirsing, E.:
{\em On the theorem of Gauss-Kusmin-L\'{e}vy and a Frobenius type theorem for function spaces},
Acta Arith. \textbf{24} (1973-74), 507--528. MR0337868

\bibitem{Wu} Wu, J.:
{\em A remark on the growth of the denominators of convergents},
Monatsh. Math., \textbf{147} (2006), 259--264. MR2215567

\bibitem{XX} Xie, J.-S.; Xu, Y.-Y.:
{\em Range-Renewal Structure of Transient Simple Random Walk},
preprint.

\end{thebibliography}



\end{document}